%

\batchmode

\def\eqdef{\buildrel \rm def \over =}

\def\ms{\medskip}

\def\phi{\varphi}

\def\su{\subseteq}
\def\a{\alpha}
\def\b{\beta}

\def\l{\lambda}
\def\k{\kappa}

\def\om{\omega}

\def\lng{\langle}
\def\rng{\rangle}

\def\cont{{2^{\aleph_0}}}



\def\cf{{\rm cf\,}}


\def\endproof#1{\hfill {\parfillskip0pt$\smiley_{\hbox{{#1}}}$\par\medbreak}}

\def\imply{\Rightarrow}
\def\iff{\Leftrightarrow}
\def\proof{\smallbreak\noindent{\sl Proof}: }

\def\Cal#1{{\cal #1}}

\newcount\itemno
\def\itm{\advance\itemno1 \item{(\number\itemno)}}
\def\ritm{\advance\itemno1 \item{)\number\itemno(}}
\def\startitm{\itemno=-1 }
\def\startaitm{\itemno=0}

\def\aitm{\advance\itemno1 
\item{(\letter\itemno)}}

\def\letter#1{\ifcase#1 \or a\or b\or c\or d\or e\or f\or g\or h\or
i\or j\or k\or l\or m\or n\or o\or p\or q\or r\or s\or t\or u\or v\or
w\or x\or y\or z\else\toomanyconditions\fi}

\newcount\secno

\magnification=1200
\baselineskip 18pt
\def\cp{{\rm cp\,}}
\def\Forb{{\rm Forb\,}}

\headline{}

\def\rest{|}

\def\wcp{{\rm wcp\, }}


\def\endproof{\hfill$\triangle$\ms}

\overfullrule0pt
\def\wcp{{\rm wcp\, }}

\def\endproof{\hfill$\triangle$}

\font\bigfont cmbx10 scaled \magstep2
\font\namefont cmbx10 scaled \magstep2
{{
\obeylines

\everypar={\hskip0cm plus 1 fil}
\bf
\parskip=0.4cm

{\bigfont On universal graphs  without cliques}
{\bigfont  or without large bipartite graphs}

\bigskip
{\rm May 1995}
\bigskip

\parskip0.1cm

\bigskip

{\namefont Menachem Kojman}
Department of Mathematics 
Carnegie-Mellon University
Pittsburgh PA, 15213
{\tt kojman@andrew.cmu.edu}

}
\vfill
\everypar{}
\rm
\leftskip2cm\rightskip2cm

ABSTRACT. For every uncountable cardinal $\l$,   suitable negations of
the Generalized Continuum Hypothesis imply:
\startaitm
\aitm For all infinite $\a$ and $\b$, there is no universal
$K_{\a,\b}$-free graphs in $\l$

\aitm For all $\a\ge 3$, there is no universal $K_\a$-free graph in $\l$ 

The instance $K_{\om,\om_1}$ for $\l=\aleph_1$ was settled by Komjath
and Pach from the principle $\diamondsuit(\omega_1)$.

\footline{\hfill}
\global\pageno0

\eject
}

\headline{\tensl\hfill M.~ Kojman: graphs\hfill}

{\bf \S0 Introduction}

The Generalized Continuum Hypothesis, GCH, is an extremely useful
assumption in infinite graph theory in general, and  in the theory of
universal graphs in particular. One consequence of the GCH is the
existence of universal graphs in all infinite powers.

In this paper {\it negations} of the GCH are used to settle a few
problem in the theory of universal graphs. Some of these problems were
treated in the past with the GCH, and some were not. 

The theory of universal graphs began with Rado's construction [R] of a
strongly universal countable graph. The research in this area has
advanced considerably since Rado's paper, mainly in studying
universality in monotone classes of graphs, or, equivalently, in
classes of the form $\Forb(\Gamma)$, all graphs omitting a set
$\Gamma$ of ``forbidden'' configurations.  A good source for the
development in this area is the survey paper [KP1] in which the
authors suggest a generalization of universality they name
``complexity'': the least number of elements in the class needed to
embed all other members in the class as induced subgraphs. The
complexity of a class is 1 exactly when a universal member exists in
the class.

In this paper omissions of infinite cliques and infinite complete
bipartite graphs are studied.
The omissions of $K_\om$ and of $K_{\om,\om}$ were studied in [DHV]
and in  [HK].  Omission of $K_\a$ for uncountable $\a$ was
treated in [KS] using the GCH.  Omissions of $K_{\a,\b}$ for $\a$
finite, $\a\le \b$, were settled in [KP] for all infinite powers $\l$
from the GCH, and the omission of $K_{\om,\om_1}$ was settled
negatively for $\l=\aleph_1$ from the  principle $\diamondsuit(\om_1)$
in the same paper.

The omission $K_\a$ for all $\a\ge 3$ and the omission of $K_{\a,\b}$
for all infinite $\a\le \b$ is settled here in all uncountable powers
from suitable negations of the GCH. This complements and extends the
results of Diestel-Halin-Vogler, Komjath-Pach and Komjath-Shelah.

\ms
\noindent
NOTATION 

Write $G_1\le G_2$ if the graph $G_1$ is isomorphic to an induced
subgraph of the graph $G_2$ and $G_1\le_w G_2$ if $G_1$ is isomorphic
to a subgraph of $G_2$. A class $\Cal G$ of graphs is {\it monotone}
if $G_1\le_w G_2\in \Cal G\imply  G_1\in \Cal G$. For a set of graphs
$\Gamma$, let $\Forb(\Gamma)$ be the class of all graphs $G$ satisfying
$H\not\le_w G$ for all $H\in \Gamma$. Every monotone class is of the
form $\Forb(\Gamma)$ for some class $\Gamma$ of graphs. 

Write $\Cal G_\l$ and $\Forb_\l(\Gamma)$ for the set of all
isomorphism types of cardinality $\l$ in $\Cal G$ and in $\Forb(\Gamma)$
respectively. Let $\cp \Cal G_\l$, the {\it complexity} of $\Cal
G_\l$, be the least cardinality of a collection $D\su \Cal G_\l$
satisfying that for all $G\in \Cal G_\l$ there exists $G'\in D$ such
that $G\le G'$. Replacing $\le$ by $\le_w$ in the last definition we
obtain $\wcp \Cal G_\l$, the {\it weak complexity} of $\Cal G_\l$. For
every class $\Cal G$ and cardinal $\l$ it holds that $\wcp \Cal
G_\l\le \cp \Cal G_\l$. The complexity $\cp \Cal G_\l$ is 1 iff there
is a strongly universal graph in $\Cal G_\l$ and similarly for $\wcp
\Cal G_\l$. 

Let $\k,\l$ be cardinals. By $\cf\k$ we denote the {\it cofinality} of
$\k$. The power set $\Cal P(\k)$ is the set of all subsets of $\k$.
By $[\k]^\l$ we denote the set of all subsets of $\k$ whose
cardinality is $\l$. Let $\cf\lng [\k]^{\l},\su \rng$, the {\it
cofinality} of the partially ordered set $\lng [\k]^\l,\su\rng$ (the
partial ordering is set inclusion), be the least cardinality of a
collection $D\su [\k]^\l$ satisfying that for all $X\in[\k]^\l$ there
exists $Y\in D$ such that $X\su Y$.

\ms
\noindent
{\bf \S1 The results}
\ms

\ms\noindent{\bf 
1,1 Definition}: Let $\theta$ be infinite. For $A\su \Cal
P(\theta)$  let the {\it incidence graph}
of $A$, denoted $\Gamma_A$, be the bipartite graph with left side
$\theta$, right side $A$ and edge relation given by $\in$, the
membership relation (a set is connected to its members by edges).

\ms\noindent{\bf 
1.2 Theorem}: 
{\sl
Suppose $\theta<\l$ are infinite cardinals and  $\cal G$
is a class of graphs that contains all  incidence
graphs $\Gamma_A$ for $A\in[\Cal P(\theta)]^\l$.
If $\cp \Cal G_\l\le 2^\theta$ then $\cf
\lng[2^\theta]^\l,\su\rng\le 2^\theta$.
}

\proof Suppose $\cal F$ is a family of graphs, each of cardinality $\l$,
such that $|\Cal F|\le 2^\theta$ and every $G\in \Cal G_\l$ is embeddable
as an induced subgraph in some member of $\Cal F$.

For every $A\in [\Cal P(\theta)]^\l$ fix an embedding $f_A:\Gamma_A\to G_A$
for some $G_A\in \Cal F$.
Given a graph $G\in \Cal F$ the number of functions $g:\theta\to G$ is
at most $\l^{\theta}\le {2^\theta}^{\theta}=2^\theta$.

For every $G\in \Cal F$ and every function $g:\theta\to G$, define

$$S(G,g)\eqdef\bigcup\{A\in [\Cal P(\theta)]^\l:G=G_A\;\&\;f_A\rest
\theta=g\}$$ 

The family $\Cal F^*=\{S(G,g):G\in \Cal F, g\in G^\theta\}$ has
cardinality $\le 2^\theta$ and covers $[\Cal P(\theta)]^\l$ because
$A\su S(G_A,f_A\rest \theta)$. Since $|\Cal P(\theta)|=2^{\theta}$,
the proof will be done once we prove that every member of $\Cal F^*$
has cardinality $\le\l$. Suppose that $x,y\in S(G,g)$ are distinct,
and let $A,B$ be such that $x\in A$, $y\in B$ and
$f_A\rest\theta=f_B\rest\theta=g$. Since $x$ and $y$ are distinct,
there is a point $z\in \theta$ such that $z\in x\iff z\notin y$. As
$f_A(z)=f_B(z)$ and both functions preserve edges and non-edges, it
follows that $\{g(z),f_A(x)\}\in E^G\iff
\{g(z),f_B(y)\}\notin E^G$. Hence $f_A(x)\not=f_B(y)$. We have shown,
then, that $f=\bigcup\{f^{-1}_A :f_A\rest \theta=g\}$ is a surjection
 from  $G$ onto $S(G,g)$, and therefore  $|S(G,g)|\le \l$.\endproof

\ms\noindent{\bf 
1.3 Remark}:  The condition $\cp \Cal G_\l\le 2^\theta$ in
1.2  can be weakened to ``there exist  $\le 
2^\theta$ many graphs, each of cardinality $\l$, not all of which
necessarily belonging to $\Cal G_\l$, such that every
member of $\Cal G_\l$ is isomorphic to an induced subgraph of at
least one of them''.

\ms\noindent{\bf 
1.4 Corollaries}: 
{\sl
Suppose $\theta$ is
infinite and $\Cal G$ is a class of graphs that contains all incidence
graphs of subsets of $\Cal P(\theta)$. Then:

\startitm
\itm If  $\cf2^\theta \le \l<2^\theta$ then $\cp \Cal G_\l>2^\theta$.

\itm If $\cf2^\theta=\theta^+$ then $\Cal G$  possesses no universal
elements in any cardinal $\l$ satisfying $\theta<\l<2^\theta$; in fact
$\cp
\Cal G_\l\ge {2^\theta}^+$.

\itm It is impossible to
 compute $\cp \Cal G_\l$ in ZFC or to prove the existence of a universal
element in $\Cal G_\l$ for all  cardinals $\l>\theta$.
}

\proof To prove (0) it is enough, by Theorem 1.2, to show that
if $\cf 2^\theta\le \l$ then $\cf\lng
[2^\theta]^\l,\su\rng>2^\theta$. This is a standard diagonalization
argument: for every list of $2^\theta$ many members of $[R]^\l$
construct in $\l$ many steps a subset of $2^\theta$ of size
$\cf2^\theta$ which is not contained in any of the members in the
list.

(1) follows from (0).

To prove (0) we recall that, by Easton's results [E], for every cardinal
$\mu$ with $\cf\mu>\theta$ it is consistent with the
axioms of set theory that GCH holds below $\theta$ and
$2^\theta=\mu$. Given any cardinal $\l>\theta$, there are infinitely
many cardinals $\mu>\l$ whose cofinality is, say, $\theta^+$. Thus
by (1) the complexity $\cp \Cal G_\l$ may assume infinitely many
different values, all larger than $\l$.\endproof

\noindent
{\sl Omitting complete subgraphs.}
We apply 1.2 to omissions of complete graphs: 

\ms\noindent{\bf 
1.5 Theorem}: 
{\sl
If $\a\ge 3$ is a cardinal then corollaries
(0)--(2) above hold for $\Forb(K_\a)$. In particular, for no
uncountable $\l$ and $\a\le \l$ can one prove from the usual axioms of
set theory the existence of a universal $K_\a$-free graph in power $\l$.
}

\proof For every $A\su \cal P(\theta)$ the incidence graph of $A$ is
$K_\a$-free for all $\a\ge 3$ and $\theta\ge\aleph_0$.\endproof

Hajnal and Komjath showed in [HK] that the complexity of
$\Forb_{\aleph_0}(K_\om)$ equals exactly $\aleph_1$ (see [KS]\S2 for a
generalization of this). This shows that $\theta<\l$ cannot be relaxed
to $\theta\le \l$ in Theorem 1.2 and in 1.4(1),(2). Komjath
and Shelah showed that from the GCH it follows that $\Forb(K_\a)$ has
a universal graph in $\l\ge \a$ iff $\cf\l<\cf\a$. Theorem
1.5 above settles the problem negatively from suitable
negations of GCH, namely for all $\l\in [\cf 2^\theta,2^\theta)$ for
some $\theta$.

One may ask whether the condition $\cf 2^\theta\le\l$ in Corollary
1.4(0) is necessary, or can be replaced by $\theta<\l$. Shelah
constructs a model of set theory in [S2] in which $\cont >\l$ for a
prescribed regular uncountable $\l$ and a universal graph (in the
class of all graphs) exists in power $\l$. This was generalized by
Mekler [M] to classes of structures including $\Forb(K_n)$ for all
$n$. Komjath and Shelah [KS] construct a model in which GCH holds up
to $\k$, $2^\k$ is large and $\cp
\Forb_\k(K_{\om_1})=\k^+$.  Since Corollaries (0)-(2) hold for the
class of all graphs, $\Forb(K_n)$ and $\Forb(K_{\om_1})$, the 
singularity assumption is needed for each of these classes.

\ms 
\noindent
{\sl Omitting complete bipartite subgraphs}

We turn now to omissions of complete bipartite graphs. Theorem
1.2 does not apply to $\Forb(K_{\a,\b})$ for infinite $\a$ and
$\b$, because incidence graphs may contain copies of $K_{\a,\b}$. But
an easy variation on the proof handles this.

Let $\theta$ be an infinite cardinal.

\ms\noindent{\bf 
1.6  Definition}:
 A family $\Cal A\su \Cal P(\theta)$ is
{\it $\theta$-almost disjoint} if $|\bigcap A|<\theta$ for every $A\in
[\Cal A]^\theta$.

The cardinal arithmetic assumption $\theta=2^{<\theta}$ implies
the existence of a $\theta$-almost disjoint $\Cal A\su \Cal P(\theta)$
of cardinality $|\Cal A|=2^\theta$. 

\ms\noindent{\bf 
1.7 Problem}: Is it true that a $\theta$-almost disjoint
family of size $2^\theta$ exists over every infinite cardinal
$\theta$?

\ms\noindent{\bf 
1.8 Fact}: If $\Cal A\su \Cal P(\theta)$ is $\theta$-almost disjoint
and $A\su \Cal A$ then the incidence graph $\Gamma_A$ is
$K_{\theta,\theta}$-free. 

\ms\noindent{\bf 
1.9 Theorem}: 
{\sl
If $\theta\le\a\le\b$  are infinite cardinals
and $2^{<\theta}=\theta$ then Corollaries (0)--(2) hold for $\Cal
G=\Forb(K_{\a,\b})$. In particular, for all uncountable $\l$ one
cannot prove in ZFC the existence of a universal $K_{\a,\b}$-free in
power $\l$ for all $\b\ge\a\ge\om$.
}

\proof It is enough to prove that Theorem 1.2 holds for all
classes $\Cal G$ that contain all $K_{\theta,\theta}$-free incidence
graphs of $A\in [\Cal P(\theta)]^\l$. Using $2^{<\theta}$ fix $\Cal
A\su \Cal P(\theta)$, $\theta$-almost disjoint of cardinality
$2^\theta$. In the proof of 1.2 consider only $A\in [\Cal
A]^\l$. For such $A$, the incidence graph $\Gamma_A$ is
$K_{\theta,\theta}$-free, and therefore belongs to $\Forb(K_{\a,\b})$.
The proof shows that $\cf \lng [A]^\l,\su\rng\le 2^\theta$. Since
$|\Cal A|\le2^\theta$, also $\cf\lng [2^\theta]^\l,\su\rng\le
2^\theta$.\endproof

By a theorem of Diestel, Halin and Vogler [DHV], for every non-empty
set $\Gamma$ so that every $G\in \Gamma$ contains an infinite path,
$\wcp_{\aleph_0} \Forb(\Gamma)>\aleph_0$. The proof generalizes readily
to give $\wcp \Forb_\l(\Gamma)\ge \l^+$. Since $K_{\om,\om}$ contains
an infinite path, putting $\Gamma=\{K_{\om,\om}\}$ we obtain from
Diestel-Halin-Vogler that there is no universal $K_{\om,\om}$-free
graph in $\l$ for all infinite cardinals $\l$. Komjath and Pach use the
principle $\diamondsuit(\om_1)$ to prove that
$\wcp\Forb_{\om_1}(K_{\om,\om_1})>\om_1$. The omission of $K_{\a,\b}$
is settled from negations of GCH for all infinite $\a\le\b$ by
1.9 above.  

\bigbreak
\noindent
{\bf Discussion}
The structure of embeddability in a pretty broad spectrum of monotone
classes is seen to be sensitive to the exponent function $\theta\mapsto
2^\theta$: there are no universal graphs in those classes in a cardinal
$\l$ belonging to an interval $[\cf 2^\theta,2^\theta)$. Shelah's
consistency results show that a tighter connection to the exponent
function, one which does not necessitate the singularity of
$2^\theta$, is not possible for the same spectrum of classes. It is
reasonable to ask if there are monotone classes of graphs in which the
complexity in power $\l$ is greater than or equal to $2^\theta$ for
some smaller $\theta$, not assuming anything about the cofinality of
$2^\theta$. The answer to this is yes. In [K] a class of graphs is
defined by forbidding a set of countable configurations, and  the
complexity at an uncountable regular $\l>\aleph_1$ is shown to be at
least $\cont$ by means of a {\it representation Theorem}, asserting the
existence of a surjective homomorphism from the relation of
embeddability over the class onto the relation of set inclusion over
all subsets of reals of cardinality $\l$.

\bigbreak
\bigbreak

\noindent
{\bf References}

\ms\noindent[DHV] R.~Diestel, R.~Halin and W.~Vogler, {\sl Some remarks on
universal graphs}, Combinatorica 5 (1985) 283--293

\ms\noindent[E] W.~B.~Easton, {\sl Powers of regular cardinals},
Annals of Mathematical logic. 1 (1970) 139--178.

\ms\noindent[HK] A.~Hajnal and P.~Komjath, {\sl Embedding graphs into colored
graphs} Trans. Amer. Math. Soc 307 (1988) 395-409

\ms\noindent[K] M.~Kojman, {\sl Representing Embeddability as set inclusion}, preprint.

\ms\noindent[KP] P.~Komjath and J.~Pach, {\sl Universal Graphs without large
bipartite subgraphs}, Mathematika 31 (1984) 282--290

\ms\noindent[KP1]  P.~Komjath and J.~Pach, {\sl Universal elements and the
complexity of certain classes of infinite graphs}, Discrete Math. 95
(1991) 255--270

\ms\noindent[KP2] P.~Komjath and J.~Pach, {\sl The complexity of a class of
infinite graphs}, Combinatorica 14(1) (1994) 121--125

\ms\noindent[KS] P.~Komjath and S.~Shelah, {\sl Universal graphs without large
cliques}, preprint

\ms\noindent[R] R.~Rado, {\sl Universal graphs and universal functions}, Acta.
Arith. 9. (1964) 331-340

\ms\noindent[S2] S.~Shelah, {\sl On universal graphs without instances of CH},
Annals of Pure. Appl. Logic 26 (1984) 75--87

\ms\noindent[S3] S.~Shelah, {\sl Universal graphs without instances of CH:
revisitied}, Israel J. Math 70 (1990) 69--81

\end